\documentstyle[a4,twoside,12pt]{amsart}

\def\myauthor{Ph\`ung H{{\accent"5E o}\kern-.28em\raise.2ex\hbox{\char'22}\kern-.20em} H{a\kern-.370em\raise.16ex\hbox{\char'47}\kern.1em}i}
\def\myAUTHOR{ PH\`UNG H{{\accent"5E O}\kern-.38em\raise.8ex\hbox{\char'22}\kern-.12em}  H{A\kern-.46em\raise.80ex\hbox{\char'47}\kern.18em}i}

\def\,{\hspace{.1ex}}
\def\lora{\longrightarrow}
\def\ot{\otimes}

\def\loma{\longmapsto}
\def\Vn{{V^{\ot n}}}

\newcommand{\bbas}{\begin{eqnarray*}}
\newcommand{\eeas}{\end{eqnarray*}}
\newcommand{\bbar}{\begin{array}}
\newcommand{\eear}{\end{array}}
\newcommand{\bbs}{\begin{displaymath}}
\newcommand{\ees}{\end{displaymath}}
\newcommand{\bb}{\begin{equation}}
\newcommand{\eqbb}{\begin{equation}}
\def\ee{\end{equation}}
\def\eqee{\end{equation}}
\def\bba{\begin{eqnarray}}
\def\eea{\end{eqnarray}}

\newtheorem{thm}{Theorem}[section]
\newtheorem{lem}[thm]{Lemma}

\newtheorem{pro}[thm]{Proposition}

\def\hi{\hat{i}}

\def\hj{\hat{j}}
\def\hm{\hat{m}}

\def\hn{\hat{n}}
\def\hk{\hat{k}}
\def\hl{\hat{l}}
\def\hp{\hat{p}}
\def\hq{\hat{q}}
\def\ha{\hat{a}}
\def\hb{\hat{b}}

\def\hs{\hat{s}}
\def\hI{\hat{I}}
\def\hJ{\hat{J}}

\def\hK{\hat{K}}
\def\hL{\hat{L}}

\def\Im{\mbox{\rm Im}\,}

\def\End{\mbox{\rm End\,}}

\def\H{{\cal H\,}}

\def\A{{\cal A\,}}

\def\P{{\cal P\,}}

\def\eee{\rule{.75ex}{1.5ex}\\[1ex]}

\def\proof{{\it Proof.\ }}

\newcommand{\tr}{\mbox{\rm tr\,}}

\renewcommand{\dim}{\mbox{\rm dim\,}}

\def\rref#1{(\ref{#1})}

\font\Fraktur=eufm10 scaled\magstep1          
   \newcommand{\fraktur}[1]{\mbox{\Fraktur #1}}  %
   \font\Fraktu=eufm7 scaled\magstep1            
   \newcommand{\fraktu}[1]{\mbox{\Fraktu #1}}    %
   \font\Frakt=eufm5 scaled\magstep1             
  \newcommand{\frakt}[1]{\mbox{\Frakt #1}}      %
   \def\frak#1{\mathchoice{\fraktur {#1}}            
                        {\fraktur {#1}}            
                        {\fraktu {#1}}             
                        {\frakt {#1}}  }           

\newcommand{\Ss}{\frak S}

\font\Bbbm=msbm10 scaled\magstep1
  \def\BBbm#1{\mbox{\Bbbm #1}}
\font\Bbbb=msbm7 scaled\magstep1
  \def\BBbb#1{\mbox{\Bbbb #1}}
\font\Bbbn=msbm5 scaled\magstep1
  \def\BBbn#1{\mbox{\Bbbn #1}}

\def\Bbb#1{\mathchoice{\BBbm {#1}}
                    {\BBbm {#1}}
                    {\BBbb {#1}}
                    {\BBbn {#1}} }
\newcommand{\bK}{{\Bbb K}}

\def\id{{\mathchoice{\mbox{\rm id}}
                    {\mbox{\rm id}}
                    {\mbox{\scriptsize\rm id}}
                    {\mbox{\tiny\rm id}} }}

\def\part{\vdash}

\def\lam{{\lambda}}

\begin{document}
\bibliographystyle{plain}
\newcommand{\sign}{{\,\rm sign\,}}
\let\integral=\int
\def\int{\displaystyle\integral}

\def\mytitle{ The Integral on Quantum Super Groups of Type $A_{r|s}$}
\def\myperemail{{\small\tt phung@@ncst.ac.vn}}
\def\myemail{{\small\tt phung@@mpim-bonn.mpg.de}}
\def\myaddress{\small\it Max-Planck Institut f\"ur Mathematik, Gottfried-Claren-Str. 26, 53225, Bonn, Germany}
\def\myperaddress{\small\it Hanoi Institute of Mathematics  P.O. Box 631, 10000 Bo Ho,  Hanoi, Vietnam}
\def\mythanks{The work was done during the author's stay at the Max-Planck-Institut f\"ur Mathematik, Bonn.}
\def\myabstract{We compute the integral on matrix quantum (super) groups of type $A_{r|s}$ and derive from it the quantum analogue of (super) HCIZ integral.}

\def\amshead{
\title{The Integral on Quantum Super Groups of Type $A_{r|s}$}
\author{ PH\`UNG H{{\accent"5E O}\kern-.38em\raise.8ex\hbox{\char'22}\kern-.12em}  H{A\kern-.46em\raise.80ex\hbox{\char'47}\kern.18em}i}
\address{\myperaddress}
\email{\myperemail}
\keywords{}
\begin{abstract}\myabstract\end{abstract}
\maketitle }

\amshead

\section*{Introduction}
A formula for the integral on the group $U(n)$ was obtained by Itzykson and Zuber \cite{iz1}. It can be given in the following form
\bba \int_{U(n)}\tr(MUNU^{-1})^k[dU]=\sum_{\lambda\in\P^n_k}\frac{d_\lambda}{r_\lambda}\Phi_\lambda(M)\Phi_\lambda(N)\label{eq1}\eea
for any hermitian matrices $M$ and $N$. $\Phi_\lambda$ is the irreducible character of $U(n)$, corresponding to $\lambda$. If $\xi_1,\xi_2,...,\xi_n$ are the eigenvalues of $M$ then $\Phi_\lambda(M)=s_\lambda(\xi_1,\xi_2,...,\xi_n)$, $s_\lambda$ are the Schur functions. $d_\lambda$ is the dimension of the irreducible module of the symmetric function $\Ss_n$, $r_\lambda$ is the dimension of the irreducible representation of $U(n)$, corresponding to partition $\lambda$. Explicitly, $d_\lambda=n!\prod_{x\in[\lambda]}h(x)^{-1}$, $r_\lambda=\prod_{x\in[\lambda]}c(x)h(x)^{-1}$, where $c(x)$ is the content, $h(x)$ is the hook-length of the box $x$ in the diagram $[\lambda]$.
This formula turns out to be a special case of a formula obtained by Harish-Chandra \cite{hc1}. The integral on the left-hand side of \rref{eq1} is therefore referred as Harish-Chandra-Itzykson-Zuber (HCIZ) integral. A super analogue of this formula was obtained by Alfaro, Medina and Urrutia \cite{amu1,amu2}, it reads
\bba\label{eq2}\int_{U(m|n)} {\rm str} (MUNU^{-1})^k[dU]=\sum_{\mu\in\P^m,\nu\in\P^n\atop{|\mu|+|\nu|=k-mn}}\frac{k!}{|\mu|!|\nu|!}\frac{d_\mu d_\nu}{r_\mu r_\nu}\Phi_\lambda(M)\Phi_\lambda(N),\eea
where $\lambda=(n^m)+\mu\cup\nu'.$
The aim of this paper is to give a quantum analogue of the super HCIZ integral. In contrary to the method of Itzykson-Zuber and its super version developed by Alfaro {\it et. al.}, which is analytic, our method is purely algebraic. Thus, fist we want to give an algebraic definition of the integral.

It is well-known that an integral on a compact group is uniquely defined by its left (or right) invariance with respect to the group action. Thus, the integral can be considered as a linear functional on the function algebra of the group with an invariance property. Notice that the algebra of polynomial function on a  compact (Lie) group is a commutative Hopf algebra. Being motivated by this fact, one can give the notion of integral on an arbitrary Hopf algebra, although, such an integral does not always exist. Thus, by definition, a left integral on a Hopf algebra $H$ is a linear form $\int$ on $H$ such that $(\id\ot \int)\Delta=\int$. Analogously, one has the notion of right integral. Generally, left and right integrals may differ. For general properties of integral on Hopf algebras, the reader is referred to \cite{sweedler}.

Originally, the integral is defined only on compact groups. The algebraically defined integral exists however on any cosemisimple Hopf algebras and on finite dimensional Hopf algebras. In particular, Hopf algebras of functions on linear reductive groups possess integrals. This is explained by the fact that each linear reductive group possesses a compact form and the corresponding integral coincides with the analytically defined one on this compact form. Having this relationship in mind, we shall find in this work the integral on the function algebras on quantum (super) groups of type $A_{r-1|s-1}$, such a quantum group is understood to be a generalization of the quantum general linear super group $GL_q(r|s)$. It is defined in terms of a Hecke symmetry of birank $(r,s)$.

Our problem of finding integrals on the funtion algebra of quantum linear supergroup $GL_q(r|s)$ is thus motivated by the HCIZ integral. On the other hand, it is an interesting problem from the point of view of Hopf algebra theory. Integrals on Hopf algebras were studied by several authors since the pioneering work of Sweedler \cite{sweedler}, see e.g. \cite{sul,lin1,doi1}. For finite dimensional Hopf algebra is is known that the integrals exist uniquely up to a scalar. However, only very few examples of infinite dimensional Hopf algrebras with integra, except cosemisimple Hopf algebras.

In  the case of quantum groups of type $A_r$, which corresponds to Hecke symmetry with birank $(r,0)$, the integral was computed in \cite{ph97b}. In this case, the function algebra is co-semisimple hence we know a priori the existence of the integral. In the general case, the function algebra is not co-semisimple. However, the formula for the integral in the former case suggests us an idea of finding an integral in the latter case.

Let $R$ be a Hecke symmetry on a finite dimensional vector super space $V$ of dimension $d$, and $H_R$ the associate Hopf algebra of function on a quantum group of type $A_{r-1|s-1}$, where $(r,s)$ is the birank of $R$, see Section \ref{sect1}. Thus, as an algebra, $H_R$ is generated by $2d^2$ generators $z_i^j,t_i^j,1\leq i,j\leq d$. Using the commutation rule on $H_R$, one can show that an element of $H_R$ can be represented as a linear combination of monomials on $z_i^j,t_k^l$ of the form $Z^J_IT_K^L:=z_{i_1}^{j_1}\cdots z_{i_p}^{j_p}t_{k_1}^{l_1}\cdots t_{k_q}^{l_q}.$ From the linearity of the integral, we see that it is sufficient to find the integral on the set of momomials of the form $Z_I^JT_L^K$. Let us denote $|I|$ the cardinal number of the sequence $I$ contents. It turns out that any integral should vanish on those monomial $Z_I^JT^L_K$ which has $|I|\neq |K|$.

The formula of the integral on $Z_I^JT^L_K$ with $|I|=|K|$ is based on an operator $P_n:\Vn\lora \Vn$, $n=|I|\neq |K|$. In Section \ref{sect2} we show that the axiom for the integral is equivalent to certain condions on $P_n$. In Section \ref{sect21} we construct $P_n$. An advantage of our method comparing to the one of Itzykson-Zuber and Alfaro et.al. is that we are able to compute the integral at every monomial function, while their method gives a formula of the integral only at certain trace-polynomials. In Section \ref{sect3} we derive a quantum analogue of HCIZ integral from our integral formula. To do this we have to introduce the notion of character of $H_R$-comodules, the latter are understood to be rational representation of the quantum group. In the last section we discuss the orthogonality relation of simple $H_R$-comodules. 

\section{Quantum groups associated to Hecke symmetries}\label{sect1}
Let $V$ be a super vector space over $\bK$, a fixed field of characteristic zero. Fix a homogeneous basis $x_1,x_2,\ldots, x_d$ of $V$. We shall denote the parity of the basis element $x_i$ by $\hi$. An even operator $R$ on $V\ot V$ can be given by a matrix $R_{ij}^{kl}$: $R(x_i\ot x_j)=x_k\ot x_lR^{kl}_{ij}.$ $R$ is an even operator means that the matrix elements $R^{ij}_{kl}$ are zero, except for those with $\hi+\hj=\hk+\hl$.
 $R$ is called {\it Hecke symmetry} if the following conditions are satisfied:
\begin{itemize}\item[i)] $R$ satisfies the Yang-Baxter equation $R_1R_2R_1=R_2R_1R_2,$ $R_1:=R\ot I$, $R_2:=I\ot R$, $I$ is the identity matrix of degree $d$.
\item[ii)] $R$ satisfies the Hecke equation $(R-q)(R+1)=0$ for some $q$ which will be assumed {\it not to be a root of unity.}
\item[iii)] There exists a matrix $P_{ij}^{kl}$ such that $P_{jn}^{im}R^{nk}_{ml}=\delta^i_l\delta^k_j.$ A matrix satisfying this condition is called {\it closed}.\end{itemize}
The matrix bialgebra $E_R$ and the matrix Hopf algebra $H_R$ are define as follows. Let $\{z^i_j,t_j^i|1\leq i,j\leq d\}$ be a set of variables, $\hat{x^i_j}=\hat{t^i_j}=\hi+\hj$. We define $E_R$ as the quotient algebra of the free non-commutative algebra, generated by $\{z^i_j|1\leq i,j\leq d\}$, by the relations
\bba\label{eq1.1} (-1)^{\hs(\hi+\hp)}R^{kl}_{ps} z^p_iz^s_j&=&(-1)^{\hl(\hq+\hk)}z^k_qz^l_nR_{ij}^{qn},\quad 1\leq i,j,k,l\leq d.\eea
Here, we use the convention of summing up by the indices that appear in both lower and upper places.
And we define the algebra $H_R$ as the quotient of the free non-commutative algebra generated by $\{z^i_j,t_j^i|1\leq i,j\leq d\}$, by the relations
\bba\label{eq1.2}(-1)^{\hs(\hi+\hp)}R^{kl}_{ps} z^p_iz^s_j&=&(-1)^{\hl(\hq+\hk)}z^k_qz^l_nR_{ij}^{qn},\quad 1\leq i,j,k,l\leq d,\\ 
\label{eq1.3}(-1)^{\hat{j}(\hj+\hk)}z^i_jt^j_k&=&(-1)^{\hat{l}(\hl+\hi)}t^i_lz^l_k=\delta^i_k,\quad 1\leq i,k\leq d.\eea

The relations in \rref{eq1.1} can be considered as the commuting rule for elements of $E_R$. For $H_R$, we have the following relations, which follow immediately form \rref{eq1.2} and \rref{eq1.3}.
\bba\label{eq1.4}&
(-1)^{\hk(\hi+\hj)}R^{pj}_{ql}z_j^it_k^l=(-1)^{\hm(\hn+\hp)}t_n^pz_q^mR_{mk}^{ni},&\\
\label{eq1.5}&
(-1)^{\hs(\hi+\hp)}R^{kl}_{ps} 
t^s_jt^p_i=(-1)^{\hl(\hq+\hk)}t^l_nt^k_q R_{ij}^{qn}.&\eea

It is easy to show that $E_R$ is a bialgebra, the coproduct on $E_R$ and $H_R$ is given by 
$$\Delta(z^i_j)=z^i_k\ot z^k_j,\quad \Delta(t^i_j)=t_j^k\ot t_k^i.$$
 $H_R$ is, in fact, a Hopf algebra. The coproduct is given by
$$\Delta(z^i_j)=z^i_k\ot z^k_j,\quad \Delta(t^i_j)=t_j^k\ot t_k^i,$$ the antipode on $H_R$ is given by 
$$S(z^i_j)=(-1)^{\hj(\hi+\hj)}t^i_j,\quad S(t_j^i)=(-1)^{\hi(\hi+\hj)}C^i_kz^k_l{C^{-1}}^l_j,$$
 where $C^i_j:=P^{il}_{jl}$. To verify the axiom of Hopf algebra for $H_R$, we need the following relation
\bba\label{eq1.6}&(-1)^{\hj(\hi+\hj)}z^l_kC^j_lt^i_j=C^i_k.&\eea

We also define $D^i_j:=P^{li}_{lj}$. The matrices $C$ and $D$ play important roles in our work, they are called {\it reflection operators.} Using the Hecke equation we can show that (cf \cite{ph97b})
\bbs CD=DC=q^{-1}-(Q^{-1}-1)\tr(C).\ees

Since we are working in the category of super vector spaces, the rule of sign effects on the coproduct. More precisely, the compatibility of product and coproduct of a super bialgebra reads
\bbs \Delta(a\ot b)=\sum_{(a)(b)}(-1)^{\ha_2\hb_1}a_1b_1\ot a_2b_2.\ees
Therefore we have
\bbs \Delta(z^{i_1}_{j_1}z^{i_2}_{j_2}\cdots z^{i_2}_{j_2})=\sign(I,K)\sign(K,J)\sign(I,J)z^{i_1}_{k_1}z^{i_2}_{k_2}\cdots z^{i_n}_{k_n}\ot z^{k_1}_{j_1}z^{k_2}_{j_2}\cdots z^{k_2}_{j_n},\ees
where $I:=(i_1,i_2,\ldots,i_n)$ and so on.
$\sign(I,J)$ is given recurently by
\bbs \sign(i,j)=1,\quad \sign(Ii,Ji)=(-1)^{\hi(|\hI|+|\hJ|)}\sign(I,J),\ees
$\hI:=(\hi_1,\hi_2,...,\hi_n)$, $|\hI|$ denotes the sum of its terms. Hence, for convenience, we denote 
$$Z^I_J={Z^{\ot n}}^I_J:=\sign(I,J)z^{i_1}_{j_1}z^{i_2}_{j_2}\cdots z^{i_2}_{j_2}.$$
Throughout this paper we shall always use this notation.
Then we have 
$$\Delta(Z^I_J)=Z^I_K\ot Z^K_L,\quad \Delta(T^I_J)=T^K_L\ot T^I_K,$$
 and $S(Z^I_J)=(-1)^{|\hJ|(|\hI|+|\hJ|)}T^{I'}_{J'}.$ Notice that 
\bbs \sign(K',L')=(-1)^{|\hat{K_1}|(\hk_n+\hj_n)}\sign(K_1',L_1'),\ees
where for $K=(k_1,k_2,...,k_n)$, $K_1:=(k_1,k_2,...,k_{n-1})$ and $K':=(k_n,k_{n-1},...,k_1)$. 

The Hecke algebra of type $A$, $\H_{n,q}$ is generated by elements $T_i,1\leq i\leq n-1$, subject to the relations
$$T_iT_{i+1}T_i=T_{i+1}T_iT_{i+1},\quad T_i^2=(q-1)T_i+q, \quad i=1,...,n-2.$$
To each element $w$ of the symmetric group $\Ss_n$ of permutations of the sets $\{1,2,...,n\}$, one can associated in a canonical way an element $T_w$ of $\H_n=\H_{n,q}$, in particular, $T_1=1, T_{(i,i+1)}=T_i$. The set $\{T_w|w\in\Ss_n\}$ form a $\bK$ basis for $\H_n$.

$R$ induces an action of the Hecke algebra $\H_n=\H_{q,n}$ on the tensor powers $\Vn$ of $V$,  $\rho_n(T_i)=R_i:=\id_V^{i-1}\ot R\ot\id_V^{n-i-1}.$ We shall therefore use the notation $R_w:=\rho(T_w)$. On the other hand,  $E_R$ coacts on  $V$ by $\delta(x_i)=x_j\ot z^j_i$. Since $E_R$ is a bialgebra, it coacts on $\Vn$ by means of the product. The double centralizer theorem \cite{ph97} asserts that these two actions are  centralizers of each other in $\End_\bK(\Vn)$ \cite{ph97}. Hence, the algebra $\End^{E_R}(\Vn)$ is a factor algebra of $\H_n$ and $(E_R^n)^*\cong \End_{\H_n}(\Vn).$ Let us denote $\End^{E_R}(\Vn)$ by $\overline{\H_n}$. Since $\H_n$ is semi-simple, provided $q$ is not a root of unity, the algebras $(E_R^n)^*$ and $\overline{\H_n}$ are semi-simple, too.

The double centralizer theorem also implies that a simple $E_R$-comodule is the image of the operator induced by a primitive idempotent of $\H_n$ and, conversely, each primitive idempotent of $\H_n$ induces an $E_R$ comodule which is either zero or simple. On the other hand, irreducible representations of $\H_n$ are parameterized by partitions of $n$. Thus, up to conjugation, primitive idempotents of $\H_n$ are parameterized by partitions of $n$, too. Note that by the semisimplicity, $\overline{\H_n}$ is also a subalgebra of $\H_n$.

The primitive idempotents $x_n:=\sum_w T_w/[n]_q!$ and  $y_n:=\sum_w (-q)^{-l(w)}T_w/[n]_{1/q}!$ induce the symmetrizer and anti-symmetrizer operators on $\Vn$. Let $S_n:=\Im\rho_n(x_n)$ and $\Lambda_n:=\Im\rho_n(y_n)$. Then one can show that $S:=\bigoplus_{n=0}^\infty S_n$ and $\Lambda:=\bigoplus_{n=0}^\infty \Lambda_n$ are algebras. They are called symmetric and exterior tensor algebras on the corresponding quantum space. 

By definition, the Poincar\'e series $P_\Lambda(t)$ of $\Lambda$ is $\sum_{n=0}^\infty t^n\dim_\bK(\Lambda_n)$. It is proved that this series is a rational function having negative roots and positive poles. Let $r$ be the number of its roots and $s$ be the number of its poles.
As a consequence, simple $E_R$-comodules are parameterized by hook-partitions from  $\Gamma^{rs}_n:=\{\lambda\part n|\lam_{r+1}\leq s\}$ \cite{ph97c}. Therefore, in the algebra $\overline{\H_n}$ primitive idempotents are parameterized by hook-partitions from $\Gamma^{r,s}_n$, too.

$(r,s)$ is called the birank of $R$. 
Our main assumption on $R$ is that
\bba\label{assumption}\tr(C)=-[s-r]_q,\eea
$C$ is the reflection operator introduced above.
It can be proved that $\tr(C)$ should have the form $-[x]_q$ for some integer $x$ in the interval $[-s,r]$ \cite{ph98a}. The above equation holds for any known Hecke symmetry. It is conjectured that it holds for all Hecke symmetry.

Simple $H_R$-comodules are much more complicated. The problem of classifying all its simple comodules is still open.

The Hopf algebra $H_R$ (resp. the bialgebra $E_R$) is called the (function algebra on) a quantum group (resp. quantum semigroup) of type $A_{r-1|s-1}$.

The following are two main examples of Hecke symmetries. The Drinfeld-Jimbo's $R$-matrix of type $A_{r-1}$ \cite{jimbo86}, for $1\leq i,j,k,l\leq r$, $p^2=q$,
\bbs {R_r}^{kl}_{ij}:=\left\{
\bbar{lll}
 q &\mbox{ if }& i=j=k=l\\
 q-1 &\mbox{ if }&k=i>j=l\\
p&\mbox{ if }&k=j\neq i=l\\
 0&\multicolumn{2}{l}{\mbox{ otherwise. }}
\eear\right.\ees

Assume that all parameters are even, the  Hopf algebra associated to $R_r$ is called the (function algebra on) standard quantum general linear group GL$_q(r)$. $R_r$ is a one-parameter deformation of the permuting operator. The super version of this operator was given by Manin \cite{manin2}. Assume that the variable $x_i$, $i\leq r$ are even and the rest $s$ variables are odd. Define, for $1\leq i,j,k,l\leq r+s$, $p^2=q$,
\bbs {R_{r|s}}^{kl}_{ij}:=\left\{
\bbar{lll}
q&\mbox{ if }& i=j=k=l, \hi=0\\
-1&\mbox{ if }& i=j=k=l, \hi=1\\
q-1 &\mbox{ if }&k=i<j=l\\
 (-1)^{\hat{i}\hat{j}}p&\mbox{ if }&k=j\neq i=l\\
 0&\multicolumn{2}{l}{\mbox{ otherwise. }}
\eear\right.\ees
$R_{r|s}$ is a deformation of the permuting operator in super symmetry.
The associated Hopf algebra is called the (function algebra on) standard quantum general linear super group GL$_q(r|s)$.

$R_r$ has the birank $(r,0)$. $R_{r|s}$ has the birank $(r,s)$.

\section{The Integral on $H_R$}\label{sect2}

Recall that by definition, a left integral on a Hopf algebra $H$ over a field $k$ is a non trivial lineal functional $\int:H\lora k$ with the invariance propety:
\bba\label{eq2.1} \int (a)=\sum_{(a)}a_{1}\ot a_2.\eea
Since we are considering super algebra, we shall also require that the integral is even, that means the value of an integral at an odd element of $H_R$ is zero. It easy to see that \rref{eq2.1} is equivalent to
\bba\label{eq2.11} \sum_{(b)}\int (aS(b_1))b_2&=\sum_{(a)}a_1\int(a_1S(b)).\eea

From the definition of $H_R$, an arbitrary element of $H_R$ can be represented
by as a linear combination of monomials in $z_i^j$ and $t_k^l$. On the other hand, using the relation on $H_R$ and the axiom (iii) of $R$, we can represent a monomial like $t^i_jz_k^l$ as a linear combination of monomials of the form $z_p^qt_r^s$, i.e., we can interchange the order of $z'$ and $t's$ in a tensor product. Namely, we have, according to \rref{eq1.5},
\bba\label{eq2.2} (-1)^{\hat l(\hat i+\hat j)}t^j_iz_k^l=(-1)^{\hat r(\hat (p+\hat q)}R_{ks}^{jp}z^q_pt^s_r P_{qi}^{rl}.\eea
Thus, by using the rule \rref{eq2.2}, we can represent any element of $H_R$ as a linear combination of monomials of the form
$$ Z_I^JT_K^L:=z_{i_1}^{j_1}\cdots z_{i_p}^{j_p}\cdots t_{k_1}^{l_1}\cdots t_{k_p}^{l_p}.$$
Therefore, it is sufficient to compute the integral at such monomials. Let us denote, for $n=|I|$,
$${I_n}_{IK}^{JL}:=(-1)^{|\hat K|(|\hat K|+|\hat J|)}\int\left(Z_I^LT_{K'}^{J'}\right)$$
where, as define in the previous section, $K'$ is the same sequence as $K$ but in the revese order. We have the following conditions on $I_n$:
\begin{itemize}\item[(i)] $I_n$ should be invariant with respect to the relation with thin $z'$s and $t'$s, given in \rref{eq1.2} and \rref{eq1.5}, respectively, that is, for all $i,j=1,2,\ldots, n-1$
$$(R_i\ot R_j)I_n=I_n(R_j\ot R_i).$$
\item[(ii)] when we contract $I_n$ with respect to the relation \rref{eq1.3}, \rref{eq1.6}, we should get $I_{n-1}$, more precisely,
\bbas \delta_{j_n}^{i_n}{I_n}_{I_1i_nK}^{J_1j_nL}&=&{I_{n-1}}^{J_1L_1}_{I_1K_1}\delta_{kn}^{l_n}\\
C_{l_n}^{k_n}{I_n}^{JL_1l_n}_{IKk_n}&=&{I_{n-1}}^{J_1j_nL_1}_{Ii_pK_1}C^{i_n}_{j_n}.\eeas

\item[(iii)] $I_n$ should respect the rule \rref{eq2.11}, which reads
$${I_n}_{IM}^{JL} Z^M_K=Z^L_N{I_n}^{JN}_{IK}.$$\end{itemize}

This condition is, in fact, sufficient for an integral on $H_R$. For, assume we have a collection of matrices $I_n$, satisfying the conditions (i-iii) above, then we can extend it linearly on the whole $H_R$. The only ambiguity that may occur is that, there may be more than one way of leading an element of $H_R$ to a linear combination of monomials of the form $Z_I^JT_K^L$. However, the Yang-Baxter equation on $R$ ensures that different ways of using rule \rref{eq2.2} give us the same result, up to relations in $z'$s and $t's$, respectively.

We thus reduced the problem to finding a family of matrices $I_n$ sastisfying conditions (i-iii). Our next claim is that $I_n$ can be found in the following way
\bba\label{eq2.3} I_n=\sum_{w\in\Ss_n}(-q)^{-l(w)}(P_nC^{\ot n}R_{w^{-1}})\ot R_w\eea
where $R_w=\rho(T_w)$ as in Section \ref{sect1}, $C$ is the reflection operator introduced in Section \ref{sect1} and $P_n$ is a certain operator on $\Vn$. More precisely, we have
\begin{lem}\label{lem2.1} Assume that the operator $P_n\in \overline{\H_n}=\End^{H_R}(\Vn)\subset\End^k(\Vn)$  is in the center of $\overline{\H_n}$ and sastisfies the condition
$$(P_{n-1}\ot \id_V)=P_n(L_n+\tr(C))$$
where $L_n$ are the Murphy operators: $L_1=0$, 
$$L_n=\sum_{i=1}^{n-1}q^{-i}R_{(n-i,n)}, \quad n\geq 2,$$
$(n-i,n)$ is the invesion that changes places of $n-i$ and $n$. Then the matrices $I_n$ given \rref{eq2.3} satisfy the conditions (i-iii).\end{lem}
\proof The conditions (i) and (iii) can be easily verified. In fact, (i) is equivalent to the equations
\bbas\lefteqn{ \sum_{w\in\Ss_n}(-q)^{-l(w)}(R_iP_nC^{\ot n}R_{w^{-1}}\ot R_w=}&&\\
&& \sum_{w\in\Ss_n}(-q)^{-l(w)}(P_nC^{\ot n}R_{w^{-1}})\ot (R_wR_i),\eeas
for $i=1,2,\ldots,n-1$.
By assumption, $P_n$ commutes with all $R_i$. On the other hand, using the Yang-Baxter equation we can also show that $C^{\ot n}$ commutes with all $R_i$. Therefore the equation above implies from the following
$$ \sum_{w\in\Ss_n}(-q)^{-l(w)}(R_iR_{w^{-1}})\ot R_w=\sum_{w\in\Ss_n}(-q)^{-l(w)}(R_{w^{-1}})\ot (R_wR_i)$$
which can be easily verified using the Hecke equation for $R$. The verification of (iii) is strightforward, it does not involve  $P_n$ and $C^{\ot n}$ but rather a direct consequence of relations \rref{eq1.2}.

The harder part is to verify (ii). Here we use the condition 
$$(P_{n-1}\ot \id_V)(L_n+\tr(C))=P_n.$$
In fact, the operator $L_n$ comes into play by the following equality.

It is known that each element $T_w$ of $\H_n$ can be expressed in the form
$T_w=T_k\cdots T_{n-1}T_{w_1}$ for some $w_1\in\Ss_{n-1}$, where $\Ss_{n-1}$ is the subgroup of $\Ss_n$, fixing $n$, i.e. $\Ss_{n-1}$ permutes only $\{1,2,...,n-1\}$. We define a linear map $\H_n\lora \H_{n-1}$ setting
$$h_n(T_w)=\left[\bbar{ll} T_k\cdots T_{n-2}T_{w_1}&\mbox{if  } k\leq n-2\\
\tr(C)T_{w_1}&\mbox{if  }k=n-1.\eear\right.$$
Then we have an identity in $\H_n$
\bba\label{eq2.4}
\sum_{w\in\Ss_n}(-q)^{-l(w)}T_{w^{-1}}\ot h_n(T_w)=\sum_{u\in\Ss_{n-1}}(-q)^{-l(w)}(L_n+\tr(C))T_{u^{-1}}\ot T_u\eea
here we identify $\H_{n-1}$ with the subalgebra of $\H_n$ generated by $T_u,u\in\Ss_{n-1}$. For the proof the reader is refered to \cite{ph97b}. In fact, we can replace $\tr(C)$ by any element of $\bK$, but the crucial point here is that
$$C^{i_n}_{j_n}{R_w}^{J_1j_n}_{I_1i_n}=h_n({R_w}^J_I).$$
Here, we specialize $h_n$ on the algebra $\overline{\H_n}$.

 When we replace $T_w$ by $R_w$ in \rref{eq2.4} we obtain  immediately the second equation in (ii). More work is needed for the first equation of (ii). Interested reader is again referred to \cite{ph97b} for detail.

Thus, we reduced the problem of finding an integral on $H_R$ to constructing operators $P_n\in\overline{\H_n}$ satisfying certain conditions. This step is motivated by the construction of the integral for $H_R$ when the birank of $R$ is $(r,0)$ in \cite{ph97b}. The essential step is to construct $P_n$. The main difficulty here is that, unlike the case of of birank $(r,0)$, the operators $L_n+\tr(C)$ are not invertible in $\overline{\H_n}$, so that we cannot follow the old way to set $P_n=(P_{n-1}\ot \id_V)(L_n+\tr(C))^{-1}.$

\section{The construction of $P_n$}\label{sect21}
 
We want to construct operators $P_n\in\overline{\H_n}$ with the property
 $$P_n(L_n+\tr(C))^{-1}=P_{n-1}\ot \id_V.$$

Originally the Murphy operators were introduced by Dipper and James \cite{dj2} to describe a full set of mutually orthogonal primitive idempotents of the algebra $\H_n$: the set
\bbs E_{t_i(\lambda)}=\prod_{{1\leq m\leq n\atop |k|\leq m-1}\atop k\neq c_{t_i(\lam)}(m)}\frac{{L_m}-[k]_q}{[c_{t_i(\lam)}(m)]_q-[k]_q}, \quad 1\leq i\leq d_\lambda,\ \lambda\in\P_n,\ees
where $c_{t_i(\lambda)}(m)$ is the content of $m$ in the standard tableau $t_i(\lambda)$, is a full set of mutually orthogonal primitive idempotents of $\H_n$. The idempotents $E_{t_i(\lam)},$ $i=1,2,...,d_\lam$ belong to the same block that corresponds to $\lambda$, their sum $\sum_{1\leq i\leq d_\lambda}E_{t_i(\lambda)}=F_\lambda$ -- the minimal central idempotent corresponding to $\lambda$.

 $L_m$  satisfy the equation
\bbs \prod_{k=-m-1}^{m+1}(L_m-[k]_q)=0.\ees
Therefore, for $1\leq m\leq n$,
\bba\label{eq32} L_mE_{t_i(\lambda)}=E_{t_i(\lambda)}L_m=c_{t_i(\lambda)}(m)E_{t_i(\lambda)}.\eea

 We define
\bbs p_\lambda:=\prod_{x\in [\lambda]\setminus[(s^r)]}\frac{q^{r-s}}{[c_\lambda(x)+r-s]_q},\ees
and
\bba\label{eq3}
P_n:=\sum_{\lambda\in\Omega_n^{r,s}\atop{1\leq i\leq d_\lambda}}p_\lambda E_{t_i(\lambda)}=\sum_{\lambda\in\Omega_n^{r,s}}p_\lambda F_\lambda.\eea

Recall that $P_n$ are defined in the algebra $\overline{\H_n}\cong\End^{E_R}(\Vn)$, which is the factor algebra of $\H_n$ by the two-sided ideal generated by minimal central idempotents corresponding to partitions from $\P_n\setminus\Gamma^{r,s}_n$. Fixing an embedding $\overline{\H_n}\hookrightarrow \overline{\H_{n+1}}$, $\overline{\H_n}\ni W\loma W\ot\id_V\in\overline{\H_{n+1}}$, we identify  $\overline{\H_n}$ with a subalgebra of $\overline{\H_{n+1}}$.
\begin{lem}\label{lem1}
The operators $P_n$ are central in $\overline{\H_n}$ and satisfy the equation (in $\overline{\H_{n+1}}$)
\bbs P_{n+1}(L_{n+1}-[s-r]_q)=P_n.\ees\end{lem}
\proof The operator $P_{n}$ is obviously central, for it is a sum of central elements.

We check the equation above. First, notice that, if $\lambda\in \Omega^{r,s}_{n+1}$ and $t_i(\lambda)$ is a standard $\lambda$-tableau, then the node of $[\lambda]$, containing $n+1$, is removable, i.e., having removed it we still have a standard tableau. The tableau $t_i(\lambda)$ is called essential if this node is not the node $(r,s)$, otherwise, it is called non-essential. A tableau $t_i(\lambda)$ is essential iff the tableau, obtained from it by removing the node containing $n+1$ is again a $\gamma$-tableau with $\gamma\in\Omega^{r,s}_n.$

Observe, that if $t_i(\lambda)$ is non-essential, then $c_{t_i(\lambda)}=s-r$, hence
\bbs E_{t_i(\lambda)}(L_{n+1}-[s-r]_q)=0,\ees
by virtue of Equation \rref{eq32}.
We reshuffle the terms of $P_{n+1}$ in groups as follows
\bbs P_{n+1}=\sum_{\gamma\in\Omega^{r,s}_n\atop{1\leq i\leq d_\gamma}}
 \sum_{\lambda\in\Omega_{n+1}^{r,s}\atop{t(\lambda)\supset t_i(\gamma)}}
p_{\lambda}E_{t(\lambda)} + 
\sum_{t(\lambda)\ \rm is\atop{ not\ eessential}}p_\lambda E_{t(\lambda)}.\ees
That is, we pick up into group for each $t_i(\gamma)$, $\gamma\in\Omega^{r,s}_n$, those standard tableaux $t(\lambda)$, that contain $t_i(\gamma)$ as a subtableau. The above observation implies that the second sum in the right-hand side of the above equation multiplied by $L_{n+1}-[s-r]_q$ vanishes. Thus, it is sufficient to prove, for a fixed $t_i(\gamma)$, $\gamma\in\Omega^{r,s}_n$,
\bbs\sum_{\lambda\in\Omega_{n+1}^{r,s}\atop{t(\lambda)\supset t_i(\gamma)}}
p_{\lambda}E_{t(\lambda)}(L_{n+1}-[s-r]_q)=q_{\gamma}E_{t_i(\gamma)}.\ees

We have $(L_{n+1}-[s-r]_q)E_{t(\lambda)}=[c_{t(\lambda)}(n+1)]_q-[s-r]_q$ and $p_\lambda(c_\lambda(x)-[s-r]_q)=p_\gamma$ whenever $[\gamma]$ is obtained from $[\lambda]$ by removing the node $x$. Since, for any two standard tableaux $t(\gamma)$ and $t(\lambda)$ with $\gamma\subset\lambda$ as above, the number $n+1$ should lie in the node $x$, for which $[\lambda]\setminus [x]=[\gamma]$, we deduce that the equation to be proved is equivalent to
\bba\label{eq5}\sum_{\lambda\in\Omega_{n+1}^{r,s}\atop{t(\lambda)\supset t_i(\gamma)}}
E_{t(\lambda)}=E_{t_i(\gamma)}.\eea
 Since $\prod_{k=-n-1}^{n+1}(L_{n+1}-[k]_q)=0$,
\bbs \sum_{m=-n-1}^{n+1}\prod_{k=-n-1,\atop{k\neq m}}^{n+1}\frac{L_{n+1}-[k]_q}{[m]_q-[k]_q}=1.\ees
Therefore
\bbs E_{t_i(\gamma)}=\sum_{m=-n-1}^{n+1}E_{t_i(\gamma)}\prod_{k=-n-1,\atop{k\neq m}}^{n+1}\frac{L_{n+1}-[k]_q}{[m]_q-[k]_q}.\ees
Remember that we are working in  the algebra $\overline{\H_{n+1}}$, in which $E_\lambda\neq 0$ iff $\lambda\in\Gamma^{r,s}_{n+1}$.
Each term on the right-hand side of the above equation is either zero or a primitive idempotent of the form $E_{t(\lambda)}$ with $t(\lambda)$ containing $t_i(\gamma)$ as a subtableau. Since the left-hand side of \rref{eq5} contains all primitive idempotents in $\overline{\H_{n+1}}$ that correspond to standard tableaux containing $t_i(\gamma)$ as a subtableau, the equation \rref{eq5} follows. Lemma \ref{lem1} is therefore proved.\eee

As a colorrary of Lemmas \ref{lem2.1} and \ref{lem1}, we have
\begin{thm}\label{thm1}The Hopf algebra $H_R$ associated to a Hecke symmetry $R$, which satisfies the condition \rref{assumption}, possesses an integral, which is uniquely determined up to a scalar multiple. Let $(r,s)$ be the birank of $R$. Then an integral can be given as follows: if $l(I)=l(K)< rs$,
$\displaystyle \int(Z_I^JT_K^L)=0,$
if $l(I)=l(K)=n\geq rs$,
\bba\label{eq33}\int(Z_I^JT_{K'}^{L'})=(-1)^{|\hK|(|\hK|+|\hL|)}\sum q^{-l(w)}(P_nC^{\ot n}R_{w^{-1}})^{L}_I{R_w}^J_{K},\eea
($K'=(k_n,k_{n-1},...,k_1)$).\end{thm}
Notice that since the operators $C,P,R_w$ are all even, the integral vanishes unless $|\hI|=|\hL|, |\hK|=|\hJ|$.

There is a symmetric bilinear form on  the Hecke algera $\H_n$, given by $(T_u,T_v)=q^{l(u)}\delta^u_{v^{-1}}$. With respect to this bilinear form, $\{ R_w,w\in \Ss_n\}$ and $\{ q^{-l(w)}R_{w^{-1}}\}$ are dual bases. Thus, if $\{ E_\lambda^{ij},\lambda\part n,1\leq i,j\leq d_\lambda\}$ is a basisof $\H_n$ with the following properties
\begin{enumerate}\item $\{ E_\lambda^{ij},\leq i,j\leq d_\lambda\}$ is a basis of the block in $\H_n$, corresponding to $\lambda$;
\item $E_\lambda^{ij}E_\mu^{kl}=\delta^\mu_\lambda\delta^j_kE_\lambda^{il}.$\end{enumerate}
then, using standard argument we can easily show
\bba\label{eq331}  \sum_{w\in\Ss_n}q^{-l(w)}R_{w^{-1}} R_w=\sum_{\lambda\part n\atop{1\leq i,j\leq d_\lambda}}\frac{1}{k_\lambda}E^{ij}_\lambda\ot E_\lambda^{ji},\eea
where $k_\lambda=(E_\lambda^{ii},E_\lambda^{ii})$, $(\cdot,\cdot)$ is the mentioned above bilinear form. $k_\lambda$ can be computed explicitly (cf. \cite{ph98})
\bbs k_\lambda=q^{n(\lambda)}\prod_{x\in[\lambda]}[h(x)]_q^{-1},\ees
where $n(\lambda)=\sum_i\lambda_i(i-1)$, $h(x)$ is the hook length of $x$ in the diagram $[\lambda]$: $h_\lambda(x)=\lambda_i+\lambda_j'-i-j+1$, where $(i,j)$ is the coordinate of $x$ in the diagram $[\lambda]$.

The formula \rref{eq33} can therefore be rewritten as follows:
\bba\label{eq34} \int(Z_I^JT_{K'}^{L'})=(-1)^{|\hK|(|\hK|+|\hL|)}\sum_{\lambda\part n\atop{1\leq i,j\leq d_\lambda}}\frac{1}{k_\lambda}{E^{ij}_\lambda}_I^L {E_\lambda^{ji}}_K^J,\eea

\section{Characters of $H_R$ and quantum analogue of super HCIZ integral}\label{sect3}
For a coribbon Hopf algebra, one associate to each finite-dimensional $H$-comodule $M$ a trace map $\Phi:\End_\bK(M)\lora H$, satisfying the following conditions, for $f,g\in\End_\bK(M)$, $h\in\End_\bK(N)$,
\bba\label{eq3.1} \Phi(f\circ g)&=&\Phi(g\circ f)\\
\label{eq3.2}\Phi(f+g)&=&\Phi(g)+\Phi(g)\\
\label{eq3.3}\Phi(f\ot g)&=&\Phi(f)\Phi(g).\eea
The reader is referred to \cite{ph98} for the explicit construction of the trace map $\Phi$. Note that the ribbon structure on $H$ is to ensure the condition \rref{eq3.3} for $\Phi.$ The character of comodule $M$ is defined to be $\Phi(\id_M)\in H$.

Our Hopf algebra $H_R$ is a corribbon Hopf algebra (cf. \cite{ph98a}). We state without proof the following lemma. The reader can also consider this lemma as the definition of the characters of $H_R$-comodules $M_\lam$ and their duals.

For each partition $\lambda\in\Gamma^{r,s}_n$, $M_\lambda$ denotes the corresponding simple $E_R$-comodule. Since the map $E_R\lora H_R$ is injective (cf. \cite{ph97b}), $M_\lambda$ is simple $H_R$-comodule, too. $M_\lambda$ is isomorphic to $\Im\rho(E_\lambda)$ of a primitive idempotent $E_\lambda$. Therefore $\Phi(M_\lambda)=\Phi(E_\lambda)$. Thus,  setting $D^i_j=r(z^i_l,S(z^l_j))$, we have 
\begin{lem}\label{lem2} The characters of $M_\lambda$ and $M_\lambda^*$, $\lambda\part n$ are given by
\bbas S_\lambda:=\Phi(M_\lambda)=q^{n(r-s+1)/2}\tr(D^{\ot n}E_\lambda Z^{\ot n}),\\
S_{-\lambda}:=\Phi(M^*_\lambda)=q^{n(r-s+1)/2}\tr(C^{\ot n}E_\lambda\overline{ T^{\ot n}}),\eeas
where $\overline{ T^{\ot n}}^I_J:=T^{I'}_{J'}.$
\end{lem}
Remember that in the definition of $Z^{\ot n}$, the signs are also inserted.

The proof of this lemma does not differ from the one for the non-super case given in \cite{ph98}.

Let $K$ be an algebra. A $K$-point of a (super) bialgebra $B$ is an algebra homomorphism $\A:B\lora K$. Let $\A$ be a $K$-point of the bialgebra $E_R$. Then the entries of $A$ commute by the same rule as the entries of $Z$. We define $S_\lambda(A):=\A(S_\lambda)$. In other words, $S_\lambda(A)$ is $S_\lambda$ computed at $Z=A$. We are now ready to formulate a formula to compute the quantum super HCIZ integral.
\begin{thm}\label{thmaqiz} Let $M$ be a $K$-point of $E_R$ and $N$ be a $K$-point of $E_R'=S(E_R)$ (i.e., the subbialgebra of $H_R$, generated by $T$). Assume that entries of $M$ and $N$ commute and that they commute with the entries of $Z$ and $T$. Then 
\bbs \int \tr(D^{\ot n}M^{\ot n}Z^{\ot n}\overline{N^{\ot n}}\overline{T^{\ot n}})=
q^{-n(r-s+1)}\sum_{\lambda\in\Omega^{r,s}_n}\frac{d_\lambda p_{\lambda}}{k_\lambda}S_\lambda(M)S_{-\lambda}(N).\ees

\end{thm}
\proof 
\bbas \lefteqn{\int\tr(D^{\ot n}M^{\ot n}Z^{\ot n} \overline{N^{\ot n}}\overline{T^{\ot n}})}&&\\
&=&\sum_{w\in\Ss_n}q^{-l(w)}\tr(P_nC^{\ot n}R_{w^{-1}}\overline{N^{\ot n}})\cdot \tr(R_wD^{\ot n}{M^{\ot n}})\\
&=& \sum_{1\leq i,j\leq d_\lambda\atop{\lambda\part n}}k_\lambda^{-1}
\tr(P_nE_\lambda^{ij}C^{\ot n}\overline{N^{\ot n}})\cdot \tr(E_\lambda^{ij}D^{\ot n}{M^{\ot n}})\\
&=&\sum_{1\leq i,j\leq d_\lambda\atop{\lambda\part n}}q^{-n(r-s+1)}k_\lambda^{-1}p_\lambda
\Phi(E_\lambda^{ij*})(N)\cdot \Phi(E_\lambda^{ji})(M)\\
&=&\sum_{\lambda\part n}q^{-n(r-s+1)}d_\lambda k_\lambda^{-1}p_\lambda
S_{-\lambda}(N)S_\lambda(M).\eeas
Since $p_\lambda=0$ whenever $\lambda\not\in\Omega_n^{r,s}$, we obtain the 
desired equality.\eee

\vskip2ex
\noindent{\bf Example.} 
Let us consider the case of standard quantum general linear super group $GL_q(r|s)$, determined in terms of the symmetry $R_{r|s}$ given Section \ref{sect1}. In this case, any diagonal matrix with commuting entries is a point of $E_R$. Thus, assume that $M$ and $N$ are diagonal matrix with entries commuting each other and with the entries of $Z$ and $T$, $M= {\rm diag}(m_1,m_2,...,m_{r+s})$, $N= {\rm diag}(n_1,n_2,...,n_{r+s})$. Then, we have
\bbs S_{(n)}=\sum_{k=0}^n h_{n-k}(qm_1,...,q^r,m_r)e_k(-q^rm_{r+1},...,-q^{r-s+1}m_{r+s}),\ees
$h_k,e_k$ are the $k^{\rm th}$ complete and elementar symmetric functions in $r$ and $s$ variables, respectively. Hence
$$S_\lambda(M)=s_\lambda(qm_1,q^2m_2,...,q^rm_r/-q^rm_{r+1},...,-q^{r-s+1}m_{r+s}),$$
$s_\lambda$ are the Hook-Schur functions in $r+s$ variables (cf. \cite{mcd2}, Ex. I.3.23).
Therefore, if $\lambda\in\Omega^{r,s}_n$, thus, $\lambda=(r^s)+\mu\cup \nu'$, $\mu\in\P^r, \nu\in\P^s$, we have [loc.cit.]
$$S_\lambda(M)=(-1)^{|\nu|}\prod_{i=1,j=1}^{r,s}(q^im_i-q^{r-j+1}m_{r+j})s_\mu(\{ q^im_i\})s_\nu(\{ q^{r+1-i}m_{r+i}\}),$$
where $s_\mu$ (resp. $s_\nu$) are the Schur functions in $r$ variables (resp. $s$ variables).
Analogously, we have
$$S_{-\lambda}(N)=(-1)^{|\nu|}\prod_{i=1,j=1}^{r,s}(q^{r-s-i+1}n_i-q^{j-s}n_{r+j})s_\mu(\{ q^{r-s+1-i}n_i\})s_\nu(\{ q^{i-s}n_{r+i}\}).$$
The quantum super HCIZ is then given by, ($n\geq rs$),
\bbas &\lefteqn{ \int\tr(C^{\ot n}N^{\ot n}Z^{\ot n}M^{\ot n}T^{\ot n})=}&\\
&\displaystyle \sum_{\mu\in\P^r,\nu\in\P^s,\atop{|\mu|+|\nu|=n-rs\atop{\lambda=(s^r)+\mu\cap \nu'}}}&q^{r|\nu|-s|\mu|}p_\mu p^-_\nu\frac{d_\lambda}{k_\lambda}
\prod_{i=1\atop{j=1}}^{r,s}(q^im_i-q^{r-j+1}m_{r+j})(q^{r-s-i+1}n_i-q^{j-s}n_{r+j})\\[1ex]
&&\times s_\mu(\{ q^im_i\})s_\nu(\{ q^{r+1-i}m_{r+i}\})s_\mu(\{ q^{r-s+1-i}n_i\})s_\nu(\{ q^{i-s}n_{r+i}\})
..\eeas
\section{The Orthogonal Relations}\label{sect4}
We are now interested in the orthogonal relations.
Let $M_\lambda, M_\mu$ be two simple comodules corresponding to partitions $\lambda$ and $\mu$ of $n$. Let $M_\lambda=\Im E_\lambda$, $M_\mu=\Im E_\mu$. Then, using Theorem \ref{thm1}, we have
\bbas (\Phi(M_\lambda),\Phi(M_\mu))&=&\int(\Phi(M_\lambda)\Phi(M_\mu^*))\\
&=&\sum q^{-l(w)}\tr(P_nC^{\ot n}E_\lambda R_{w^{-1}}E_\mu R_w)\\
&=& \delta_\lambda^\mu\delta^\lambda_{\Omega^{r,s}}\frac{p_\lambda}{k_\lambda}\tr(C^{\ot n}E_\lambda)\\&=&0.\eeas
Here $\delta^\lambda_{\Omega^{r,s}}$ indicates, whether $\lambda$ belongs to $\Omega^{r,s}$, it is zero if $\lambda\not\in\Omega^{r,s}$ and 1 otherwise. On the other hand, for $\lambda$ form $\Omega^{r,s}$, $\tr(C^{\ot n}E_\lambda)=0$, for it is the quantum rank of $M_\lambda$.

We see that the scalar product above cannot be used to define the orthogonal relations. Following \cite{amu1} we compute the integral $\displaystyle\int(Z_\lambda Z_{-\mu})$, where $Z_\lambda$ is a coefficient matrix of the simple comodule $M_\lambda$, i.e., $Z_\lambda$ is the multiplicative matrix corresponding to  certain basis of $M_\lambda$.

Let us fix a primitive idempotent $E_\lambda$, $\lambda\part n$, and set $M_\lambda=\Im\rho(E_\lambda)$. To choose a basis of $M_\lambda$ we proceed as follows. First notice that there exists a one-one correspondence between the set of bases of $M_\lambda$ and the set of maximal operators $P^a_b$ in $\End_\bK(M_\lambda$ with property $P^a_bP^c_d=\delta^c_bP^a_d.$ Further, since $E_\lambda$ is an idempotent, $\End_\bK(M_\lambda)=E_\lambda\End_\bK(\Vn)E_\lambda$. Therefore, instead of fixing a basis of $M_\lambda$ we can equivalently choose a maximum set of operators $P^a_b=E_\lambda Q^a_bE_\lambda$, for some $Q^a_b\in\End(\Vn)$, satisfying $P^a_bP^c_d=\delta^c_bP^a_d.$

Let $C_\lambda$ be the restriction of $C^{\ot n}$ on $M_\lambda$, which is an invariant space of $C^{\ot n}$) and ${C_\lambda}^a_b$ be the matrix element of $C_\lambda$ with respect to the basis give by $P_b^a$, $1\leq a,b\leq m_\lambda$, $m_\lambda$ is the dimension of $M_\lambda$. Then we have
$$ {C_\lambda}_b^a=\tr(C^{\ot n}P^a_b).$$
Analogously, the multiplicative matrix, corresponding the basis determined by $P_b^a$, is given by ${Z_\lambda}^a_n=\tr(P^a_bZ^{\ot n})$. For the dual comodule we have ${Z_{-\lambda}}^a_b=\tr(P^a_b\overline{T^{\ot n}})$. Now we can find the integral
\bbas \int({Z_\lambda}^a_b{Z_{-\lambda}}^c_d)&=&\int \tr(P^a_bZ^{\ot n})\tr(P^c_d\overline{T^{\ot n}})\\
&=& k^{-1}_\lambda p_\lambda\tr(C^{\ot n}E^{ij}_\lambda P^a_bE^{ji}_\lambda P^c_d)\\
&=& k^{-1}_\lambda p_\lambda\delta^c_b\tr(C^{\ot n}P^a_d)\\
&=&k^{-1}_\lambda p_\lambda\delta^c_b{C_\lambda}^a_d,\eeas
for we can choose $E^{ij}_\lambda$ such that $E_\lambda=E^{ii}_\lambda$ for some $i$. Thus, we have proved
\begin{pro}\label{pro41} Let $M_\lambda$ be the simple comodule of $H_R$, corresponding to partition $\lambda$ and $e_1,e_2,\cdots,$ be its basis. Let $Z_\lambda$ be the corresponding multiplicative matrix and $Z_{-\lambda}$ be the multiplicative matrix corresponding to the dual basis on ${M_\lambda}^*$. Let $C_\lambda$ be the restriction of the operator $C^{\ot n}$ on $M_\lambda$ and ${C_\lambda}_b^a$ be its matrix elements with respect to the basis above. The we have the following orthogonal-type relations:
\bba\label{eq41} 
 \int({Z_\lambda}^a_b{Z_{-\lambda}}^c_d)=k^{-1}_\lambda p_\lambda\delta^c_b{C_\lambda}^a_d,\eea\end{pro}

\begin{center}\bf Acknowledgement\end{center}

\mythanks

\end{document}